\documentclass{article}
\usepackage{amssymb}
\usepackage{amsthm}
\usepackage{amsmath}
\newtheorem{Defn}{Definition}[section]
\newtheorem{Lem}[Defn]{Lemma}

\newtheorem{Thm}[Defn]{Theorem}

\begin{document}
\renewcommand{\thefootnote}{}
\title{On Sharp Local Turns of Planar Polynomials}
\author{Ruixiang Zhang}
\date{}
\maketitle
\textbf{Abstract}\ We show that for a real polynomial of degree $n$ in two variables $x$ and $y$, any local ``sharp turn'' must have its ``size'' $\gtrsim e^{-Cn^{2}}$. We also show that there is indeed an example that has a sharp turn of size $\lesssim e^{-Cn}$. This gives a quite satisfactory answer to a problem raised by Guth. The formulation of the problem was inspired by applications of the polynomial method in the study of Kakeya conjecture.

\section{Introduction} The polynomial method has led to a lot of progress in Kakeya-type problems. Perhaps the most famous one was the positive answer to the finite field Kakeya conjecture by Dvir \cite{Dvir}. However when one works on $\mathbb{R}$, people don't know whether the polynomial method will lead to any progress of the Kakeya conjecture, though a ``multilinear'' version was shown by the polynomial method thanks to the work of Guth \cite{Guthendpoint}. The difficulty with the real numbers comes from some loss of the algebraic structure due to replacing lines by tubes (for this one can look at the last part of this lecture note by Guth \cite{Guthclass}). Thus it is conceivable that some new ingredients have to show up in this setting if we want to prove the Kakeya conjecture in $\mathbb{R}^n$ by the polynomial method. And it is very likely to be some geometric properties of an arbitrary real polynomial.

One possible direction is to estimate $\frac{|\bigcup 10T_j|}{|\bigcup T_j|}$, where $T_j$'s are thin tubes of different directions. It is natural to ask how a polynomial can behave in $\bigcup 10T_j$ if it bisects each small cube in $\bigcup T_j$. This is a ``global'' question and one can first look at its ``local'' version: how sharp can a polynomial turn near one point \cite{Guthclass}?

Very little was known even for this ``local toy model'' \cite{Guthclass} and the purpose of this paper is to make a much better understanding about it than before. Namely, we prove the following quantitative theorem about local ``sharp turns'' of real polynomials in two variables. This theorem partially answers a problem very similar to, and morally the same as a problem raised in \cite{Guthclass}.

\begin{Thm}\label{mainthm}
Suppose a polynomial $f(x, y)$ with real coefficients and degree $n$ satisfies that for some $\varepsilon >0$,
\begin{equation}
f(x, y) \left\{ \begin{array}{ll}
>0 & \text{ if } \varepsilon \leq x \leq 1, -1 \leq y \leq 1\\
>0 & \text{ if } \varepsilon \leq y \leq 1, -1 \leq x \leq 1\\
<0 & \text{ if } -1 \leq x, y \leq -\varepsilon
\end{array} \right.
\end{equation}

Then we must have $\varepsilon \gtrsim e^{-Cn^{2}}$. Moreover, for any $n$ there exists an example with $\varepsilon \lesssim e^{-Cn}$ for some universal constant $C>0$.
\end{Thm}

If we visualize this theorem, we see it gives estimates about the best a planar real polynomial can do to approximate a right angle (locally). And that's where the ``sharp turn'' in the title comes from. Clearly, it is a ``toy'' version of the bisecting problem we mentioned before.

This theorem is somewhat surprising because one may expect that no examples essentially better than the simple example $(x+1)^{2n} + (y+1)^{2n} = 1$ can be made. However Theorem \ref{mainthm} says that one can do much better (exponentially). This is basically a bad news in our $\frac{|\bigcup 10T_j|}{|\bigcup T_j|}$ setting since it means that even if a polynomial bisects each small cube in a long tube $T$, it is possible that it doesn't bisect a positive portion of small cubes in $2T\backslash T$, which in turn means we can almost achieve nothing by a ``local'' analysis in this setting.

However, it is still conceivable that a polynomial is unlikely to have a lot of sharp turns. This means that it is  possible that some ``global'' results much better than Theorem \ref{mainthm} still holds. If that is the case, it may be helpful to our  $\frac{|\bigcup 10T_j|}{|\bigcup T_j|}$ problem. Unfortunately, I don't even know how to formulate a global version conjecture.

In my point of view, the method we are going to use to prove Theorem \ref{mainthm} is also interesting in its own right. Our proof has two parts. To find a lower bound of $\varepsilon$, we notice that any monomial $w$ in two variables cannot make the set $\{ (x, y) : w(x, y) > 0\}$ look like $\{(x, y) : x > 0 \text{ or } y>0\}$ near the origin. This is trivial by analysing whether $w$ is odd or even in the two variables. We shall make this property work for us and deduce that a linear combination of these monomials could not make its positive locus near the region above, either. To construct an example for the second part of Theorem \ref{mainthm}, we superpose two polynomials, one is even in $y$ and large (in absolute value) when $\varepsilon \leq x \leq 1$ and the other is odd in $y$ and large (in absolute value) when $-1 \leq x \leq -\varepsilon$.

In our discussion, all polynomials will have real coefficients. We use $C$ to denote a positive constant that may vary even in a single chain of inequalities.

\section*{Acknowledgements}

The author was supported by mathematics department of Princeton University. He would like to thank Ben Yang for bringing this interesting problem to his attention, and Yuan Cao for helpful discussions which made him notice the advantage of considering the inner inscribed ball. He would also like to thank the helpful referee whose advice made the exposition of this paper improved.

\section{A lower bound of $\varepsilon$}

In this section we establish the lower bound estimate in Theorem \ref{mainthm}. Namely, we will prove
\begin{Lem}\label{lowerbd}
Under the assumptions of Theorem \ref{mainthm}, $\varepsilon \gtrsim e^{-Cn^{2}}$.
\end{Lem}

We rotate the $xy$ -  plane by $\frac{3}{4} \pi$ counterclockwisely around the origin and easily sees that Lemma \ref{lowerbd} holds if we can prove the following lemma. For our convenience, we still denote the rotated polynomial by $f$.

\begin{Lem}\label{annuluslem}
Suppose a polynomial $f(x, y)$ with real coefficients and degree $n$ satisfies that for some $\varepsilon >0$,
\begin{equation}
f(x, y) \left\{ \begin{array}{ll}
<0 & \text{ if } \frac{1}{2} \leq r \leq 1 \text{ and } \frac{\sqrt{2}}{2} + \varepsilon \leq \cos \theta \leq 1\\
>0 & \text{ if } \frac{1}{2} \leq r \leq 1 \text{ and } -1 \leq   \cos \theta \leq  \frac{\sqrt{2}}{2} - \varepsilon
\end{array} \right.
\end{equation}
where $x+\text{i}y = re^{\text{i} \theta}$. Then we have $\varepsilon \gtrsim e^{-Cn^{2}}$.
\end{Lem}

Obviously, under the assumptions of Lemma \ref{annuluslem}, we may assume that $f(x, y)$ is even in $y$. Otherwise we consider $f(x, y) + f(x, -y)$. Using polar coordinates, we see that $f(x, y)$ is a polynomial in $r, \cos \theta$ and $\sin \theta$. Since $f$ is even in $y$, all powers of $\sin \theta$ are even. Thus $f (x, y) = \phi (r, cos \theta)$ where $\phi$ is a polynomial of degree $\leq n$ both in $r$ and in $\cos \theta$.

We observe one more property of this polynomial: Assume $\phi (r, \cos \theta) = \sum_{j = 0}^n r^j \phi_j (\cos \theta)$. By the definition of $\phi$ it is immediate that $\phi_j$ is odd (in $\cos \theta$) when $j$ is odd, and is even when $j$ is even.

Denote $s = \cos \theta$. Then as a consequence of this discussion, we showed that to prove Lemma \ref{annuluslem}, it suffices to prove the following lemma.

\begin{Lem}\label{finalversion}
Suppose a polynomial $\phi (r, s) = \sum_{j = 0}^n r^j \phi_j (s)$ satisfies that the $\phi_j$ 's have degree $\leq n$, each $\phi_j$ is either even or odd, and that  for some $\varepsilon >0$,
\begin{equation}
\phi(r, s) \left\{ \begin{array}{ll}
<0 & \text{ if } \frac{1}{2} \leq r \leq 1 \text{ and } \frac{\sqrt{2}}{2} + \varepsilon \leq s \leq 1\\
>0 & \text{ if } \frac{1}{2} \leq r \leq 1 \text{ and } -1 \leq  s \leq  \frac{\sqrt{2}}{2} - \varepsilon
\end{array} \right.
\end{equation}
Then we have $\varepsilon \gtrsim e^{-Cn^{2}}$.
\end{Lem}

\begin{proof}
First we sketch the philosophy of the proof. Note that for each summand $\phi_j (s)$ , we have $|\phi_j (s)| = |\phi_j (-s)|$. We are going to use this and the conditions to establish an estimate of the kind $|\Phi (s)| \lesssim_n |\Phi (-s)|$ where $\Phi$ will be a (nonvanishing) positive linear combination of $\phi_j (s), 0 \leq j \leq n$ and $|s| \leq 1, |s-\frac{\sqrt{2}}{2}| \geq \varepsilon$. Also our $\Phi $ will be $<0$ if $\frac{\sqrt{2}}{2} + \varepsilon \leq s \leq 1$ and $>0$ if $-1 \leq  s \leq  \frac{\sqrt{2}}{2} - \varepsilon$. But when $\varepsilon \rightarrow 0$ this two properties will become impossible to hold simultaneously because heuristically, the ``growth rate'' of $\Phi$ will be different near $-\frac{\sqrt{2}}{2}$ and $\frac{\sqrt{2}}{2}$ due to a simple sign change analysis. Here everything can be made effective and our lower bound will drop out.

Now we come to details. To avoid repetition, we introduce a notation. If a polynomial $\psi (s)$ satisfies that
\begin{equation}
\psi(s) \left\{ \begin{array}{ll}
<0 & \text{ if }  \frac{\sqrt{2}}{2} + \varepsilon \leq s \leq 1\\
>0 & \text{ if }  -1 \leq  s \leq  \frac{\sqrt{2}}{2} - \varepsilon
\end{array} \right.,
\end{equation}
then we say that $\psi$ has the ``\emph{sharp-change-at-$\frac{\sqrt{2}}{2}$}'' property.

Due to the assumptions, every $F_{r_0} (s) = \phi (r_0 , s)  = \sum_{j = 0}^n r_0^j \phi_j (s)$ has this property when $\frac{1}{2} \leq r_0 \leq 1$. In particular, $\sum_{j=0}^n a_j \phi_j (s)$ has the sharp-change-at-$\frac{\sqrt{2}}{2}$ property if $(a_0, a_1, \ldots, a_n) = (1, \frac{n+k}{2n}, \ldots, (\frac{n+k}{2n})^n)$, $\forall 0 \leq k \leq n$.

Furthermore, we notice a crucial principle that if two polynomials both have the sharp-change-at-$\frac{\sqrt{2}}{2}$ property, an arbitrary positive linear combination of them will also have this property. Thus, any point $(a_0, a_1, \ldots, a_n)$ inside the $n+1$ dimensional simplex $\Delta_n$ with vertices $(1, \frac{n+k}{2n}, \ldots, (\frac{n+k}{2n})^n) (0 \leq k \leq n)$ and the origin will make the polynomial $\sum_{j=0}^n a_j \phi_j (s)$ have the sharp-change-at-$\frac{\sqrt{2}}{2}$ property.

Now we choose $(a_0, a_1, \ldots, a_n)$ to be the center $(b_0, b_1, \ldots, b_n)$ of the inner inscribed ball $B_n$ of $\Delta_n$. The advantage of this choice is that if we denote $\Phi (s)  = \sum_{j=0}^n b_j \phi_j (s)$, then $\Phi (s)$ as well as $\Phi (s) + \lambda \phi_j (s)$ will have the sharp-change-at-$\frac{\sqrt{2}}{2}$ property whenever $0 \leq j \leq n$ and $|\lambda| \leq \frac{R_n}{2}$. Here $R_n > 0$ is the radius of $B_n$.

Therefore, $\Phi (s)$ and $\Phi (s) + \lambda \phi_j (s)$ will be $>0$ when $-1 \leq  s \leq  \frac{\sqrt{2}}{2} - \varepsilon$ and $<0$ when $\frac{\sqrt{2}}{2} + \varepsilon \leq s \leq 1$. This will be impossible if $|\phi_j (s)| > \frac{2}{R_n} |\Phi (s)|$ for any $j$ and any $|s| \leq 1, |s -\frac{\sqrt{2}}{2} | \geq \varepsilon$. Therefore we deduce
\begin{equation}
|\phi_j (s)| \leq \frac{2}{R_n} |\Phi (s)|, 0\leq j \leq n, |s| \leq 1, |s -\frac{\sqrt{2}}{2} | \geq \varepsilon.
\end{equation}

But we notice that $|\phi_j (s)| = |\phi_j (-s)|$. And it is trivial that in the expression $\Phi (s)  = \sum_{j=0}^n b_j \phi_j (s)$ each $|b_j| \leq 1$. So we have
\begin{equation}\label{likeevenorodd1}
|\Phi (-s)| \leq \frac{2(n+1)}{R_n} |\Phi (s)|, 0\leq j \leq n, |s| \leq 1, |s -\frac{\sqrt{2}}{2} | \geq \varepsilon.
\end{equation}

We postpone the detailed calculation of $R_n$ to the end of the proof where we shall see it is a small number of magnitude $O(e^{-Cn})$. Now let's simply denote $c_n = \frac{2(n+1)}{R_n}$ and rewrite (\ref{likeevenorodd1}) as
\begin{equation}\label{likeevenorodd2}
|\Phi (-s)| \leq c_n |\Phi (s)|, 0\leq j \leq n, |s| \leq 1, |s -\frac{\sqrt{2}}{2} | \geq \varepsilon.
\end{equation}

Without loss of generality, we may assume $\varepsilon < 100000^{-n}$. Take a positive parameter $K > 100$ such that $K^{n+1} \varepsilon < 0.1$. We draw $n+1$ disjoint pairs of annuli $A_h^-$ and $A_h^+$ ($0 \leq h \leq n$) on the complex plane such that: Each $A_h^-$ is the annulus centered at $-\frac{\sqrt{2}}{2}$ and has inner radius $K^h \varepsilon$ and outer radius  $K^{h+1} \varepsilon$. While each $A_h^+$ is the annulus centered at $\frac{\sqrt{2}}{2}$ and has inner radius $K^h \varepsilon$ and outer radius  $K^{h+1} \varepsilon$. i.e. each $A_h^+$ is the reflection image of the corresponding $A_h^-$ about $y$-axis. Notice that by the prescribed sign of $\Phi$, $\Phi$ has to be a nonzero polynomial. By pigeonhole, there exists an $h = t$ such that there are no zeros of $\Phi$ inside the annuli $A_t^-$ or $A_t^+$.

Now assume $\Phi (s) = M \prod_{j=1}^n (s - s_j)$. By assumption that $\Phi$ has no roots inside $A_t^+$, we deduce that each $s_j$ satisfies that either $|s_j - \frac{\sqrt{2}}{2}| \leq K^t \varepsilon$ or $|s_j - \frac{\sqrt{2}}{2}| \geq K^{t+1} \varepsilon$. Denote $\Gamma^+ = \{ 1 \leq j \leq n : |s_j - \frac{\sqrt{2}}{2}| \leq K^t \varepsilon\}$ and $\gamma^+ = \# \Gamma^+$. We then have
\begin{equation}\label{positiveratio}
\left| \frac{\Phi (\frac{\sqrt{2}}{2} + \frac{K^{t+1}\varepsilon}{2})}{\Phi (\frac{\sqrt{2}}{2} + 2K^t \varepsilon)} \right| = \prod_{j=1}^n \left| \frac{\frac{\sqrt{2}}{2} - s_j + \frac{K^{t+1}\varepsilon}{2}}{\frac{\sqrt{2}}{2} - s_j + 2K^t \varepsilon}\right|.
\end{equation}

Look at each multiplicand above. We easily deduce that for $j \in \Gamma^+$ the corresponding  multiplicand is between $\frac{K}{10}$ and $10K$, while for $j \notin \Gamma^+$ the corresponding  multiplicand is between $\frac{1}{10}$ and $10$. Therefore we have
\begin{equation}\label{positiveratio}
e^{-Cn}K^{\gamma^+} \lesssim \left| \frac{\Phi (\frac{\sqrt{2}}{2} + \frac{K^{t+1}\varepsilon}{2})}{\Phi (\frac{\sqrt{2}}{2} + 2K^t \varepsilon)} \right| \lesssim e^{Cn} K^{\gamma^+}
\end{equation}

Similarly, suppose $\Gamma^- = \{ 1 \leq j \leq n : |s_j + \frac{\sqrt{2}}{2}| \leq K^t \varepsilon\}$ and $\gamma^- = \# \Gamma^-$. We deduce
\begin{equation}\label{nagetiveratio}
e^{-Cn}K^{\gamma^-} \lesssim \left| \frac{\Phi (-\frac{\sqrt{2}}{2} - \frac{K^{t+1}\varepsilon}{2})}{\Phi (-\frac{\sqrt{2}}{2} - 2K^t \varepsilon)} \right| \lesssim e^{Cn} K^{\gamma^-}
\end{equation}

Note that $\Phi$ satisfies sharp-change-at-$\frac{\sqrt{2}}{2}$ property. Thus it is $>0$ when $-1 \leq  s \leq  \frac{\sqrt{2}}{2} - \varepsilon$ and $<0$ when $\frac{\sqrt{2}}{2} + \varepsilon \leq s \leq 1$. By this we deduce that there are an odd number of (real) roots in the close interval $[\frac{\sqrt{2}}{2} - K^t \varepsilon, \frac{\sqrt{2}}{2} + K^t \varepsilon ]$. Since the imaginary roots appear in pairs, we conclude that $\gamma^+$ is odd. Similarly, $\gamma^-$ is even.

Hence, $\gamma^+ \neq \gamma^-$. By this and (\ref{positiveratio}), (\ref{nagetiveratio}) we conclude that
\begin{eqnarray}\label{notsameincreasespeed}
&\left|\frac{\Phi (\frac{\sqrt{2}}{2} + \frac{K^{t+1}\varepsilon}{2})}{\Phi (-\frac{\sqrt{2}}{2} - \frac{K^{t+1}\varepsilon}{2})}\frac{\Phi (-\frac{\sqrt{2}}{2} - 2K^t \varepsilon)}{\Phi (\frac{\sqrt{2}}{2} + 2K^t \varepsilon)}\right| \gtrsim e^{-Cn}K \nonumber\\
\text{ or } &\left|\frac{\Phi (\frac{\sqrt{2}}{2} + \frac{K^{t+1}\varepsilon}{2})}{\Phi (-\frac{\sqrt{2}}{2} - \frac{K^{t+1}\varepsilon}{2})}\frac{\Phi (-\frac{\sqrt{2}}{2} - 2K^t \varepsilon)}{\Phi (\frac{\sqrt{2}}{2} + 2K^t \varepsilon)}\right| \lesssim e^{Cn}\frac{1}{K}
\end{eqnarray}

However by (\ref{likeevenorodd2}), we have
\begin{equation}\label{likeincreasingspeed}
c_n^{-2} \lesssim \left|\frac{\Phi (\frac{\sqrt{2}}{2} + \frac{K^{t+1}\varepsilon}{2})}{\Phi (-\frac{\sqrt{2}}{2} - \frac{K^{t+1}\varepsilon}{2})}\frac{\Phi (-\frac{\sqrt{2}}{2} - 2K^t \varepsilon)}{\Phi (\frac{\sqrt{2}}{2} + 2K^t \varepsilon)}\right| \lesssim c_n^2.
\end{equation}

Take $K = (\frac{0.01}{\varepsilon})^{\frac{1}{n+1}}$ and compare (\ref{notsameincreasespeed}) and (\ref{likeincreasingspeed}), we finally deduce that
\begin{equation}\label{epsilonlower}
\varepsilon \gtrsim e^{-Cn^2}c_n^{-(2n+2)} \gtrsim e^{-Cn^2}R_n^{2n+2}.
\end{equation}

Finally, we estimate $R_n$ which is the radius of the inner inscribed ball of the $n+1$ dimensional simplex $\Delta_n$. Geometrically,
\begin{equation}
R_n = \frac{(n+1)\text{Vol} (\Delta_n) }{\text{Area} (\Delta_n)}.
\end{equation}
where $\text{Area} (\Delta_n)$ is the whole area of the surface of $\Delta_n$.

Recall that $\Delta_n$ has $n+2$ vertices that are $(0, 0, \ldots, 0)$ and $(1, \frac{n+k}{2n}, \ldots, (\frac{n+k}{2n})^n)$ $(0 \leq k \leq n)$ in $\mathbb{R}^{n+1}$. Using the formula of the volume of a simplex, we have
\begin{equation}\label{volofdeltan}
\text{Vol} (\Delta_n) =
\frac{1}{(n+1)!}\left\| \begin{array}{cccc}
1 & \frac{n}{2n} & \ldots & (\frac{n}{2n})^n\\
1 & \frac{n+1}{2n} & \ldots & (\frac{n+1}{2n})^n\\
 &  & \ldots & \\
1 & \frac{2n}{2n} & \ldots & (\frac{2n}{2n})^n
\end{array} \right\|= \frac{1}{(n+1)!} \prod_{n\leq j<k\leq 2n} (\frac{k-j}{2n}).
\end{equation}

Now we estimate $\text{Area} (\Delta_n)$. Denote $P_k = (1 , \frac{n+k}{2n} , \ldots,  (\frac{n+k}{2n})^n) (0 \leq k \leq n)$. Then the area of the surface $OP_0 P_1 \cdots \widehat{P_l} \cdots P_n$ ($0 \leq l \leq n$) is no more than the sum of the area of  its projection to each coordinate hyperplane. That is to say
\begin{equation}\label{areaofkthface1}
\text{Area} (OP_0 P_1 \cdots \widehat{P_l} \cdots P_n) \leq \sum_{h=0}^n \frac{1}{n!}\left\| \begin{array}{cccccc}
1 & \frac{n}{2n} & \ldots & \widehat{(\frac{n}{2n})^h} & \ldots & (\frac{n}{2n})^n\\
1 & \frac{n+1}{2n} & \ldots &  \widehat{(\frac{n}{2n})^h} & \ldots & (\frac{n+1}{2n})^n\\
 &  & & \ldots & & \\
 \widehat{1} & \widehat{\frac{n+l}{2n}} & \ldots &  \widehat{(\frac{n+l}{2n})^h} & \ldots & \widehat{(\frac{n+l}{2n})^n}\\
  &  & & \ldots & & \\
1 & \frac{2n}{2n} & \ldots & \widehat{(\frac{n}{2n})^h} & \ldots & (\frac{2n}{2n})^n
\end{array} \right\|
\end{equation}

Using Newton's identities we can easily transfer RHS of (\ref{areaofkthface1}) to some multiples of standard Vandermonde determinants. We deduce
\begin{equation}\label{areaofkthface2}
\text{Area} (OP_0 P_1 \cdots \widehat{P_l} \cdots P_n) \leq \sum_{h=0}^n \frac{\sigma_{h, \hat{l}}}{n!}\prod_{0\leq j< k \leq n, j, k \neq l} (\frac{k-j}{2n}).
\end{equation}

Here $\sigma_{h, \hat{l}}$ is the $h$-th symmetric sum of $\frac{n}{2n}, \frac{n+1}{2n}, \ldots, \widehat{\frac{n+l}{2n}}, \ldots, \frac{2n}{2n}$. Clearly $\sigma_{h, \hat{l}} < {n\choose h}$. Thus
\begin{equation}\label{areaofkthfacefinal}
\text{Area} (OP_0 P_1 \cdots \widehat{P_l} \cdots P_n) \leq \frac{2^n}{n!}\prod_{0\leq j< k \leq n, j, k \neq l} (\frac{k-j}{2n}).
\end{equation}

Now there is one more surface, namely $P_0 P_1 \cdots P_n$. It is parallel to a coordinate hyperplane and has area
\begin{eqnarray}\label{areaofkthface0}
&\text{Area} (P_0 P_1 \cdots P_n) \nonumber\\
= & \frac{1}{n!}\left\| \begin{array}{cccc}
\frac{n+1}{2n} - \frac{n}{2n} & (\frac{n+1}{2n})^2 - (\frac{n}{2n})^2 & \ldots & (\frac{n+1}{2n})^n - (\frac{n}{2n})^n\\
\frac{n+2}{2n} - \frac{n}{2n} & (\frac{n+2}{2n})^2 - (\frac{n}{2n})^2 & \ldots & (\frac{n+2}{2n})^n - (\frac{n}{2n})^n\\
 &  & \ldots & \\
 \frac{2n}{2n}  - \frac{n}{2n} & (\frac{2n}{2n})^2 - (\frac{n}{2n})^2 & \ldots & (\frac{2n}{2n})^n - (\frac{n}{2n})^n
\end{array} \right\| \nonumber\\
= & \frac{1}{n!} \cdot \frac{n!}{(2n)^n} \left\| \begin{array}{cccc}
1 & \frac{n+1}{2n} + \frac{n}{2n} & \ldots & (\frac{n+1}{2n})^{n-1} + (\frac{n+1}{2n})^{n-2}\cdot \frac{n}{2n} + \cdots + (\frac{n}{2n})^{n-1}\\
1 & \frac{n+2}{2n} + \frac{n}{2n} & \ldots & (\frac{n+2}{2n})^{n-1} + (\frac{n+2}{2n})^{n-2}\cdot \frac{n}{2n} + \cdots + (\frac{n}{2n})^{n-1}\\
 &  & \ldots & \\
1 & \frac{2n}{2n} + \frac{n}{2n} & \ldots & (\frac{2n}{2n})^{n-1} + (\frac{2n}{2n})^{n-2} \cdot \frac{n}{2n} + \cdots + (\frac{n}{2n})^{n-1}
\end{array} \right\| \nonumber\\
= & \frac{1}{n!} \cdot \frac{n!}{(2n)^n} \left\| \begin{array}{ccccc}
1 & \frac{n+1}{2n} & (\frac{n+1}{2n})^2 & \ldots & (\frac{n+1}{2n})^{n-1}\\
1 & \frac{n+2}{2n} & (\frac{n+2}{2n})^2 & \ldots & (\frac{n+2}{2n})^{n-1}\\
 &  & & \ldots & \\
1 & \frac{2n}{2n} & (\frac{2n}{2n})^2 & \ldots & (\frac{2n}{2n})^{n-1}
\end{array} \right\| \nonumber\\
= & \frac{1}{n!} \prod_{0\leq j< k \leq n} (\frac{k-j}{2n}).
\end{eqnarray}

Now combine (\ref{volofdeltan}), (\ref{areaofkthfacefinal}) , (\ref{areaofkthface0}) and Stirling's asymptotic formula, we deduce
\begin{equation}\label{finalRn}
R_n = \frac{(n+1)\text{Vol} (\Delta_n) }{\text{Area} (\Delta_n)} \gtrsim e^{-Cn}.
\end{equation}

(\ref{epsilonlower}) and (\ref{finalRn}) together imply the lemma.
\end{proof}

As discussed in the beginning, Lemma \ref{finalversion} implies the lower bound estimate , namely Lemma \ref{lowerbd}, needed  in Theorem \ref{mainthm}.

\section{An upper bound of $\varepsilon$}

We construct a polynomial to achieve the upper bound in Theorem \ref{mainthm}.

\begin{Lem}\label{upperbd}
Under the assumptions of Theorem \ref{mainthm}, it is possible that $\varepsilon \lesssim e^{-Cn}$.
\end{Lem}

\begin{proof}
Without loss of generality, we assume $n> 100$.

Our plan is to construct some positive constant $C_0$ and polynomials $p(y), q(y)$ and $g(x)$ of degree $\leq \frac{n-1}{2}$ which are in one variable (respectively) and satisfy the following:
\begin{equation}\label{cond1}
p(y), q(y) \text{ are even functions of } y.
\end{equation}
\begin{equation}\label{cond2}
p(y), q(y) > 0, y \in \mathbb{R}.
\end{equation}
\begin{equation}\label{cond3}
1 \leq \frac{p(y)}{yq(y)} \leq 1.1 , e^{-C_0 n } \leq y\leq 1.
\end{equation}
\begin{equation}\label{cond4}
p(y) \geq yq(y) , 0 \leq y \leq 1.
\end{equation}
\begin{equation}\label{cond5}
g(x) \geq 0, x \in \mathbb{R}.
\end{equation}
\begin{equation}\label{cond6}
g(x) > 2 g (-x), e^{-C_0 n } \leq x\leq 1.
\end{equation}

Once (\ref{cond1})-(\ref{cond6}) are satisfied, it is easy to verify that the polynomial $p(y) g(x) + yq(y) g(-x)$ satisfies all the conditions of Theorem \ref{mainthm} when $\varepsilon = e^{-C_0 n}$.

Therefore it remains to construct $p(y), q(y), g(x)$ such that (\ref{cond1}) - (\ref{cond6}) hold. We first construct $p$ and $q$ inductively. Take a parameter  $a \in (0, 1)$ (to be fixed later) and construct a series of even polynomials $(p_n (y), q_n (y))$ inductively as follows:

First we set $p_0 (y) = q_0 (y) = 1$. Hence $\frac{p_0}{y q_0}$ maps the interval $[a, 1]$ into $[1, \frac{1}{a}]$ and $p_0 (y) \geq yq_0 (y)$ when $0 \leq y \leq 1$.

Suppose we already have $p_m, q_m$ even, such that $p_m (y) \geq yq_m (y)$ and that $p_m (y), q_m (y)  >0$ when $0 \leq  y \leq 1$ and $\frac{p_m}{y q_m}$ maps $[a, 1]$ into $[1, K_m]$. We construct $p_{m+1} (y) = p_m^2 (y) + K_m y^2 q_m^2 (y)$, $q_{m+1} (y) = 2 \sqrt{K_m} p_m (y) q_m (y)$. Then it is immediate to check that $p_{m+1}$ and $q_{m+1}$  are even and that $p_{m+1} (y) , q_{m+1} (y)  >0$ when $0 \leq  y \leq 1$. Also, $\frac{p_{m+1} (y)}{y q_{m+1} (y)} = \frac{1}{2} (\frac{p_m (y)}{\sqrt{K_m} y q_m (y)} + \frac{\sqrt{K_m} y q_m (y)}{p_m (y)})$ for $y>0$. Thus $p_{m+1} (y) \geq yq_{m+1} (y)$ when $0 \leq  y \leq 1$. And by induction hypothesis, $\frac{p_{m+1}}{y q_{m+1}}$ maps $[a, 1]$ into $[1, K_{m+1}] = [1, \frac{1}{2} (\sqrt{K_m} + \frac{1}{\sqrt{K_m}})]$.

Note that $\max \{\deg p_{m+1} , \deg{q_{m+1}}\} \leq 2 \max \{\deg p_m , \deg{q_m}\} +2 $. Thus if we choose $m_0$ to be the largest $m$ such that $\max \{\deg p_{m} , \deg{q_{m}}\} \leq \frac{n-1}{2}$, we have $m_0 \gtrsim \frac{\log n}{\log 2} -C$. We take $a = 1.1^{-2^{m_0}}$. Note that $K_{m+1} < \sqrt{K_m}$, we deduce that $\frac{p_{m_0}}{y q_{m_0}}$ maps $[a, 1] = [1.1^{-2^{m_0}}, 1]$ into $[1, K_{m_0}] \subseteq [1, (\frac{1}{a})^{\frac{1}{2^{m_0}}}] \subseteq [1, 1.1]$. Since $a = 1.1^{-2^{m_0}}$, $m_0 \gtrsim \frac{\log n}{\log 2} -C$, we deduce $a \leq e^{- C_1 n}$ where $C_1$ is some positive constant. If we choose $p = p_{m_0}$, $q = q_{m_0}$ and $C_0 < C_1$ , then (\ref{cond1}), (\ref{cond2}), (\ref{cond3}) and (\ref{cond4}) are satisfied.

It suffices to concoct some $g$. We let $g(x) = \prod_{m=1}^{\lfloor\frac{n}{10}\rfloor } (x+ \frac{1}{2^m})^2$. Then its degree $\leq \frac{n-1}{2}$ and (\ref{cond5}) is automatically satisfied. Also, (\ref{cond6}) holds if $g(-x) =0$. We claim that when $g(-x) \neq 0$ and $\frac{1}{2^{\lfloor\frac{n}{10}\rfloor }} < x <1$, (\ref{cond6}) holds. Indeed, in this case
\begin{equation}
\frac{g(x)}{g(-x)} = \prod_{m=1}^{\lfloor\frac{n}{10}\rfloor } (\frac{x+ \frac{1}{2^m}}{x-  \frac{1}{2^m}})^2.
\end{equation}

All the multiplicands above are $>1$ for positive $x$. Moreover, if $\frac{1}{2^j} < x < \frac{1}{2^{j-1}} (1 \leq j \leq\lfloor\frac{n}{10}\rfloor)$, then the multiplicand $(\frac{x+ \frac{1}{2^j}}{x-  \frac{1}{2^j}})^2 > 4$. Thus (\ref{cond6}) holds for all $\frac{1}{2^{\lfloor\frac{n}{10}\rfloor }} < x <1$. This means if we choose $C_0$ small enough such that $e^{-C_0 n} > \frac{1}{2^{\lfloor\frac{n}{10}\rfloor}}$, then $g(x)$ will satisfy the conditions (\ref{cond5}) and (\ref{cond6}).

Thus by taking $C_0$ sufficiently small we can find suitable $p, q$ and $g$ satisfying  (\ref{cond1})-(\ref{cond6}) above and $p(y) g(x) + yq(y) g(-x)$ will be an example which verifies the lemma.

\end{proof}

Department of Mathematics, Princeton University, Princeton, NJ 08540

ruixiang@math.princeton.edu

\end{document}